\documentclass[10pt]{amsart}
\usepackage{amsfonts,amssymb,amsmath}
 
\newtheorem{thm}{Theorem}[section] 
\newtheorem{lem}[thm]{Lemma} 
\newtheorem{prop}[thm]{Proposition}
\newtheorem{cor}[thm]{Corollary}
\newtheorem{conj}[thm]{Conjecture}
\newtheorem*{thmA}{Theorem A}
\newtheorem*{thmB}{Theorem B}
\newtheorem*{thmC}{Theorem C}

\theoremstyle{remark}
\newtheorem{remark}[thm]{Remark}

\def\threep{\mathcal{I}_{g,3}^*}
\def\clos{\overline{\mathcal{I}}_{g,3}^*}
\def\norm{ {}' \clos}

\numberwithin{equation}{section}

\begin{document}

\title{Symmetric roots and admissible pairing}
\author{Robin de Jong}


\begin{abstract}  Using the discriminant modular form and the Noether formula it
is possible to  write the admissible self-intersection of the relative dualising
sheaf of a semistable hyperelliptic curve over a number field or function field
as a sum, over all places, of a certain adelic invariant $\chi$.  We provide a
simple geometric interpretation for this invariant $\chi$, based on the
arithmetic of  symmetric roots. We propose the conjecture that the invariant $\chi$ coincides with the invariant $\varphi$ introduced in a recent paper by S.-W. Zhang.  
\end{abstract}

\maketitle
\thispagestyle{empty}

\section{Introduction}

Let $X$ be a hyperelliptic curve of genus $g \geq 2$ over a field $K$ which is
either a number field or the function field of a curve over a field.  Assume
that $X$ has semistable reduction over $K$. We study, 
for each place $v$ of $K$, a real-valued invariant $\chi(X_v)$ of $X
\otimes K_v$, with the following two properties:
\begin{itemize}
\item[(i)] $\chi(X_v)=0$ if $v$ is non-archimedean and $X$ has good reduction at
$v$;
\item[(ii)] for the admissible self-intersection of the relative dualising sheaf
$(\omega,\omega)_a$ of $X$ the formula
\[   (\omega, \omega)_a = \frac{2g-2}{2g+1} \sum_v \chi(X_v) \log Nv   \]
holds. Here $v$ runs over the places of $K$, and the $Nv$ are usual local factors 
related to the product formula for $K$.
\end{itemize}
In the function field context, the invariant $\chi$ already appears in work
of A. Moriwaki \cite{mo} and K. Yamaki \cite{ya2}, albeit in disguise. 
It follows from their work that $\chi(X_v)$
is strictly positive if $X$ has non-smooth 
reduction at $v$. In fact they prove a precise lower bound for $\chi(X_v)$ 
in terms of the
geometry of the special fiber at $v$. If $X$ has a non-isotrivial model, 
by property (ii) this yields as a corollary 
an effective proof of the
Bogomolov conjecture for $X$, i.e. the strict positivity of $(\omega,\omega)_a$. 

The precise definition of $\chi$ is given in Section~\ref{chi}. It involves the discriminant
modular form of weight $8g+4$, suitably normalised, the $\varepsilon$-invariant of
S.-W. Zhang \cite{zh2}, and the $\delta$-invariant appearing in the Noether formula for
smooth projective curves over $K$.

Our purpose is to give a geometric interpretation of the invariant $\chi$. 
Fix, for each place $v$ of $K$, an algebraic closure
$\bar{K}_v$ of $K_v$. Endow each $\bar{K}_v$ with a standard absolute value $|\cdot|_v$ (see
Section \ref{mainresult}). Then we prove:
\begin{thmA} Let $\langle \omega,\omega \rangle$ be the Deligne self-pairing of
the dualising sheaf $\omega$ of $X$ on $\mathrm{Spec}(K)$. There exists a
canonical section $q$ of $(2g+1)\langle \omega,\omega \rangle$ on
$\mathrm{Spec}(K)$, obtained by pullback from the moduli stack of smooth
hyperelliptic curves of genus $g$ over $\mathbb{Z}$, such that the equality 
\[ -\log |q|_a = (2g-2)\chi(X_v)   \]
holds for each place $v$ of $K$. Here $|\cdot|_a$ is Zhang's admissible norm 
on $(2g+1)\langle \omega,\omega \rangle$ at $v$.
\end{thmA}
The construction of $q$ yields the following simple formula
for $\chi$. 
\begin{thmB} Assume that $K$ does not have characteristic~$2$. 
Let $v$ be a place of $K$, and let $w_1,\ldots,w_{2g+2}$ on 
$X \otimes \bar{K}_v$ be the 
Weierstrass points of $X \otimes \bar{K}_v$.
Then for each $i=1,\ldots,2g+2$ the formula 
\[ \chi(X_v) = -2g \left( \log|2|_v + \sum_{k \neq i} (w_i,w_k)_a \right) \]
holds, where $(,)_a$
is Zhang's admissible pairing on $\mathrm{Div}(X \otimes \bar{K}_v)$.
\end{thmB}
In a recent paper \cite{zh2} S.-W. Zhang introduced, for any
smooth projective geometrically connected curve $X$ of genus at least~$2$ over $K$,
an invariant $\varphi(X_v)$ for each $X\otimes K_v$ such that property (i)
holds for $\varphi$, and property (ii) holds for $\varphi$ if $X$ is
hyperelliptic. We propose the conjecture that $\varphi$ and $\chi$ are
equal for all hyperelliptic curves over $K$ and all places of $K$. 
This conjecture turns out to be true in the
case $g=2$. As we will explain below, this gives a new proof of the Bogomolov
conjecture for curves of genus~$2$ over number fields.

The main tools in this paper are moduli of (pointed) 
stable hyperelliptic curves and the arithmetic of symmetric roots of $X$.
These symmetric roots were extensively studied by J. Gu\`ardia \cite{gu} in
the context of an effective Torelli theorem for hyperelliptic period matrices. 
They are defined as follows. 
Let $\kappa$ be any field of characteristic not equal to~$2$ 
and let $X$ be a hyperelliptic curve of genus $g \geq 2$ over
$\kappa$. Fix a separable algebraic closure $\bar{\kappa}$ of $\kappa$. 
Let $w_1,\ldots,w_{2g+2}$ on $X
\otimes \bar{\kappa}$ be the Weierstrass points of $X\otimes \bar{\kappa}$, 
and let $h \colon X \to
\mathbb{P}^1_\kappa$ be the quotient under the hyperelliptic involution of $X$.
Fix a coordinate $x$ on $\mathbb{P}^1_{\bar{\kappa}}$, and suppose that $w_i$ gets mapped to
$a_i$ on $\mathbb{P}^1_{\bar{\kappa}}$. Take a pair $(w_i,w_j)$ of distinct 
Weierstrass points,
and let $\tau$ be an automorphism of $\mathbb{P}^1_{\bar{\kappa}}$ such that
$\tau(a_i)=0$, $\tau(a_j)=\infty$, and $\prod_{k \neq i,j} \tau(a_k)=1$. The
$\tau(a_k)$ for $k\neq i,j$ are finite, non-zero, 
and well-defined up to a common scalar from $\mu_{2g}$, the set of
$2g$-th roots of unity in $\bar{\kappa}$. The resulting subset of $S_{2g}
\setminus \mathbb{A}_{\bar{\kappa}}^{2g} / \mu_{2g}$ is denoted by 
$\{ \ell_{ijk} \}_{k \neq i,j}$ and is called the set of \emph{symmetric roots} on $X$ and
$(w_i,w_j)$. This is clearly an invariant of $X$ and the pair
$(w_i,w_j)$ over $\bar{\kappa}$. It is easily
checked that the formula:
\[ \ell_{ijk} = \frac{ a_i - a_k}{a_j-a_k} \sqrt[2g]{ \prod_{r \neq i,j}
\frac{(a_j - a_r)}{ (a_i-a_r) } } \]
holds for each $k \neq i,j$. This formula of course has to be interpreted appropriately
if one of the $a_i,a_j,a_k$ equals infinity. 
For each given $k \neq i,j$, the  
element $\ell_{ijk}^{2g}$ of $\bar{\kappa}$ lies in
the field of definition inside $\bar{\kappa}$ of the triple $(w_i,w_j,w_k)$.

Our main result Theorem \ref{main} gives a description of a power of
$\ell_{ijk}^{2g}$ as a rational function on the moduli stack of hyperelliptic
curves with three marked Weierstrass points.
As a corollary of this result we obtain the following remarkable formula,
expressing the norm of a symmetric root as a special value of 
Zhang's admissible pairing on divisors.
\begin{thmC}  Assume that $K$ does not have characteristic~$2$. 
Let $v$ be a place of $K$, and let $w_i,w_j,w_k$ be three distinct 
Weierstrass points on $X\otimes \bar{K}_v$. 
Let $\ell_{ijk}$ in $\bar{K}_v^\times/\mu_{2g}$ be the symmetric root on
the triple $(w_i,w_j,w_k)$. Then the formula 
\[ (w_i - w_j, w_k)_a = -\frac{1}{2} \log | \ell_{ijk} |_v \]
holds.
\end{thmC}
We mention that for an archimedean place $v$ this result says that
\[ \frac{G_v(w_i,w_k)}{G_v(w_j,w_k)} = \sqrt{|\ell_{ijk}|_v} \, , \]
where $G_v$ is the Arakelov-Green's function on the compact Riemann surface
$X(\bar{K}_v)$. This formula is remarkable since it says that a special value
of some transcendental function on $X(\bar{K}_v)$ is algebraic. We believe
that this fact merits further attention. 

The main result and Theorem C are proven in Section \ref{mainresult}. After
introducing the $\chi$-invariant in Section \ref{chi} we give  in Section
\ref{intrinsic} our more intrinsic approach to $\chi$ and prove Theorems A and
B. In Section \ref{phi} we compare $\chi$ with Zhang's $\varphi$-invariant.
Throughout, the reader is assumed to be familiar with the theory of admissible
pairing on curves as in \cite{zh1}. All schemes and algebraic stacks in this
paper are assumed locally noetherian and separated.

\section{Admissible pairing and relative dualising sheaf} \label{admpairing}

We begin with a useful description of $\langle \omega ,\omega \rangle$ for
families of semistable hyperelliptic curves. The contents of this section 
straightforwardly generalise those of \cite[Section~1]{bmmb} which treats the
genus~$2$ case.

We start with a few definitions. Let $S$ be a (locally noetherian, separated) 
scheme. A proper flat family $\pi \colon X \to S$ of curves of genus $g \geq 2$ 
is called a \emph{smooth hyperelliptic curve} over $S$ if $\pi$ is smooth and admits an
involution $\sigma \in \mathrm{Aut}_S(X)$ such that $\sigma$ restricts to a
hyperelliptic involution in each geometric fiber of $\pi$. If $\pi$ is a smooth
hyperelliptic curve, the involution $\sigma$ is uniquely determined. We call
$\pi \colon X \to S$ a \emph{generically smooth semistable hyperelliptic curve}
if $\pi$ is semistable, and there exist an open dense subscheme $U$ of $S$ and
an involution $ \sigma \in \mathrm{Aut}_S(X)$ such that $X_U$ together with the
restriction of $\sigma$ to $X_U$ is a smooth hyperelliptic curve over $U$.
Again, if $\pi$ is a generically smooth semistable hyperelliptic curve, the
involution $\sigma$ is unique; we call $\sigma$ the \emph{hyperelliptic 
involution} of $X$ over $S$.

Let $\pi \colon X \to S$ be a generically smooth semistable hyperelliptic curve
of genus $g \geq 2$. Let $\omega$ be the relative dualising sheaf of $\pi$ and
let $W$ be a $\sigma$-invariant section of $\pi$ with image in the smooth locus
of $\pi$. The image of $W$ in $X$ induces a relative Cartier divisor on $X$
which we also denote by $W$. We make the convention that whenever a Cartier 
divisor on a scheme is given, the associated line bundle will be denoted by the
same symbol. Moreover we use additive notation for the tensor product of line
bundles. 

Assume for the moment that $\pi$ is smooth. By \cite[Lemma 6.2]{dj1} there
exists a unique isomorphism:
\[ \omega \xrightarrow{\cong} (2g-2)W - (2g-1) \pi^* \langle W,W \rangle \, , \]
compatible with base change, such that pullback along $W$ induces the
adjunction isomorphism:
\[ \langle W, \omega \rangle \xrightarrow{\cong} - \langle W,W \rangle \]
on $S$. Here $\langle \cdot, \cdot \rangle$ denotes
Deligne pairing of line bundles on $X$. Now drop the condition that $\pi$ is
smooth. By the above isomorphisms we have a canonical non-zero 
rational section $s$ of the line bundle:
\[ \omega - (2g-2)W + (2g-1) \pi^* \langle W,W \rangle \]
on $X$. Denote by $V$ its divisor on $X$; then $V$ is disjoint from the smooth
fibers of $\pi$, and $W^*V = \langle W,V \rangle $ is canonically trivial on
$S$. 

Next let $(W_i,W_k)$ be a pair of $\sigma$-invariant sections of $\pi$ with image in the
smooth locus of $\pi$. Denote by $V_i,V_k$ the associated Cartier divisors
supported in the non-smooth fibers of~$\pi$. We define a line bundle $Q_{ik}$ on
$S$ associated to $(W_i,W_k)$ as follows:
\[ Q_{ik} = -4g(g-1) \langle W_i,W_k \rangle - \langle W_i,V_k \rangle
- \langle V_i,W_k \rangle + \langle V_i,V_k \rangle \, . \]
Note that $Q_{ik}$ has a canonical non-zero rational section $q_{ik}$. 
We have a canonical symmetry isomorphism 
\begin{equation} \label{symmetry} 
Q_{ik} \xrightarrow{\cong} Q_{ki} 
\end{equation}
sending $q_{ik}$ to $q_{ki}$. If $S$ is the spectrum of a discrete valuation
ring we have the following functoriality of $(Q_{ik},q_{ik})$ in passing from
$\pi \colon X \to S$ to a minimal desingularisation $\rho \colon X' \to X$ of
$X$ over $S$: the sections $W_i,W_k$ lift to $\sigma$-invariant sections of
$X'$, and one obtains the relative dualising sheaf of $X'$ over $S$ as the
pullback of $\omega$ along $\rho$. It follows that the 
$V_i,V_k$ of $X'$ over $S$ are obtained by pullback as well, and so the
formation of $Q_{ik}$ and its canonical rational section $q_{ik}$ are compatible
with the passage from $X$ to $X'$.

Assume that $S$ is an integral scheme.
\begin{prop} \label{Qik}  There exists a
canonical isomorphism:
\[ \varphi_{ik} \colon \langle \omega, \omega \rangle \xrightarrow{\cong} Q_{ik}
\]
of line bundles on $S$, compatible with any dominant base change. Let $K$ be
either a complete discrete valuation field or $\mathbb{R}$ or $\mathbb{C}$ and let
$\bar{K}$ be an algebraic closure of $K$. If $S =
\mathrm{Spec}(\bar{K})$ then $\varphi_{ik}$ is an isometry for the
admissible metrics on both $\langle \omega, \omega \rangle$ and $Q_{ik}$.
\end{prop}
\begin{proof} By construction of $V_i$ we have a canonical isomorphism:
\begin{equation} \label{first}
 \omega \xrightarrow{\cong} (2g-2)W_i + V_i -(2g-1) \pi^* \langle W_i,W_i
\rangle 
\end{equation}
on $X$. Hence we have:
\begin{equation} \label{second}
 (2g-2)(W_i-W_k) \xrightarrow{\cong} V_k - V_i + (2g-1) \pi^* \langle W_i,W_i
\rangle - (2g-1) \pi^* \langle W_k,W_k \rangle \, , 
\end{equation}
canonically. Using pullback along $W_i$ and $W_k$ we find canonical
isomorphisms:
\begin{align*}
(2g-2) \langle W_i-W_k, W_i-W_k \rangle & \xrightarrow{\cong}  \langle V_k - V_i,
W_i - W_k \rangle \\
 & \xrightarrow{\cong}  \langle W_i,V_k \rangle + \langle V_i,W_k \rangle \, .
\end{align*}
Hence we find:
\begin{equation} \label{sum}
 (2g-2) ( \langle W_i,W_i \rangle + \langle W_k,W_k \rangle )
\xrightarrow{\cong} \langle W_i,V_k \rangle + \langle V_i,W_k \rangle + 4(g-1)
\langle W_i,W_k \rangle \, . 
\end{equation}
Also from (\ref{first}) one obtains:
\[ -(2g-1) \langle W_k, \pi^* \langle W_i,W_i \rangle \rangle
\xrightarrow{\cong} - \langle W_k,W_k \rangle - (2g-2)\langle W_i,W_k \rangle -
\langle V_i, W_k \rangle \]
by using the adjunction isomorphism:
\begin{equation} \label{adj}
 \langle W_k, \omega \rangle \xrightarrow{\cong} - \langle W_k,W_k \rangle \, , 
\end{equation}
and likewise:
\[ -(2g-1) \langle W_i, \pi^* \langle W_k,W_k \rangle \rangle
\xrightarrow{\cong} - \langle W_i,W_i \rangle - (2g-2)\langle W_i,W_k \rangle -
\langle W_i, V_k \rangle \, . \]
Adding these two isomorphisms, multiplying by $2g-2$, and using (\ref{sum}) 
we obtain:
\begin{align*}  
-(2g-2)(2g-1)( \langle W_i, \pi^* \langle W_k,W_k \rangle \rangle + \langle
W_k, \pi^* \langle W_i,W_i \rangle \rangle ) \xrightarrow{\cong} \\
-(2g-1) ( \langle W_i,V_k \rangle + \langle V_i,W_k \rangle ) - 4(g-1)(2g-1)
\langle W_i,W_k \rangle \, . 
\end{align*}
Now note that:
\begin{align*}
\langle \omega , \omega \rangle 
& \xrightarrow{\cong}  (2g-2)^2 \langle W_i,W_k \rangle + (2g-2)( \langle
W_i,V_k \rangle + \langle V_i,W_k \rangle ) + \langle V_i,V_k \rangle \\
&  -(2g-2)(2g-1)( \langle W_i, \pi^* \langle W_k,W_k \rangle \rangle + 
\langle W_k, \pi^* \langle W_i,W_i \rangle \rangle ) \, . 
\end{align*}
Plugging in the previous result gives the required isomorphism. 

To see
that $\varphi_{ik}$ is an isometry for the admissible metrics on both sides if
$S = \mathrm{Spec}(\bar{K})$ with $K$ a complete discretely valued field or $K
=\mathbb{R}$ or $\mathbb{C}$ it suffices to verify that both
(\ref{first}) and the adjunction isomorphism (\ref{adj}) are isometries for the
admissible metrics. But that the adjunction isomorphism (\ref{adj}) is an isometry for the admissible metrics
is precisely in \cite{zh1}, sections~2.7 (archimedean case) 
and~4.1 (non-archimedean case). That (\ref{first}) is an isometry for the admissible
metrics on both sides can be seen as follows. The curvature form (see \cite{zh1},
section~2.5) of both left and right hand side of (\ref{first}) is equal to $(2g-2)\mu$, where
$\mu$ is the admissible metric on the reduction graph of $X$ (if $K$ is a complete 
discretely valued field) or the Arakelov $(1,1)$-form on $X(\bar{K})$ (if
$K=\mathbb{R}$ or $\mathbb{C}$). This implies (see once again \cite{zh1}, section~2.5)
that the quotient of the admissible metric on $\omega$ and the metric put on $\omega$
via the isomorphism (\ref{first}) is constant on $X\otimes \bar{K}$. By restricting the
isomorphism (\ref{first}) to $W_i$ on $X\otimes \bar{K}$ we find the adjunction isomorphism
$ \langle W_i, \omega \rangle \xrightarrow{\cong} - \langle W_i,W_i \rangle $. 
As this is an isometry, so is (\ref{first}).
\end{proof}

\section{Main result} \label{mainresult}

In this section we prove our main result and derive Theorem C from it. Let~$S$ 
be an integral scheme and let $\pi \colon X \to S$ be a generically smooth
semistable hyperelliptic curve of genus $g \geq 2$ over $S$.

Assume that a triple $(W_i,W_j,W_k)$ of $\sigma$-invariant sections of $\pi$ 
is given with image in the smooth locus of $\pi$. Assume as well that the generic characteristic
of $S$ is not equal to~$2$. We view the element $\ell_{ijk}^{2g}$, 
defined fiber by fiber along the non-empty open subscheme of $S$ where 
$\pi$ is smooth and the residue characteristic is not~$2$, as a 
rational section of the structure sheaf $\mathcal{O}_S$ of $S$. 
On the other hand, from Proposition \ref{Qik} we obtain a canonical isomorphism:
\[ \psi_{ijk} = \varphi_{ik} \otimes \varphi_{jk}^{-1} \colon -Q_{ik} + Q_{jk}
\xrightarrow{\cong} \mathcal{O}_S \, , \]
compatible with dominant base change, and isometric for the admissible metrics
on both sides. This yields a rational section of $\mathcal{O}_S$ by taking the image of
$q_{ik}^{-1} \otimes q_{jk}$ under $\psi_{ijk}$. 
Our main result is that the image of $q_{ik}^{-1} \otimes q_{jk}$
under $\psi_{ijk}$ is essentially a power of 
$\ell_{ijk}^{2g}$.
\begin{thm} \label{main}  $\psi_{ijk}$ maps the rational section
$q_{ik}^{-1} \otimes q_{jk}$ of $-Q_{ik} + Q_{jk}$ to the rational section
$(-\ell_{ijk}^{2g})^{g-1}$ of $\mathcal{O}_S$.
\end{thm}
A first step in the proof is the following result.
\begin{prop} \label{firststep}
Assume that $S$ is the spectrum of a discrete valuation ring $R$,
that $X$ is regular, and that $\pi$ is smooth if the residue characteristic of
$R$ is equal to~$2$. Then the formula:
\[ -\nu(q_{ik}) + \nu(q_{jk}) = 2g(g-1) \nu(\ell_{ijk}) \]
holds, where $\nu(\cdot)$ denotes order of vanishing 
along the closed point of $S$.
\end{prop}
\begin{proof} We recall that:
\[ Q_{ik} = -4g(g-1) \langle W_i,W_k \rangle - \langle W_i,V_k \rangle
- \langle V_i,W_k \rangle + \langle V_i,V_k \rangle \, . \]
As $\langle (2g-2)W_i - \omega + V_i, V_k \rangle$ is canonically trivial we
have a canonical isomorphism:
\[ Q_{ik} \xrightarrow{\cong} -4g(g-1) \langle W_i,W_k \rangle - (2g-1) \langle
W_i,V_k \rangle - \langle V_i,W_k \rangle + \langle V_k, \omega \rangle \]
and hence:
\[ -\nu(q_{ik}) + \nu(q_{jk}) = 4g(g-1) (W_i-W_j,W_k) + (2g-1)(W_i-W_j,V_k) +
(V_i-V_j,W_k) \, , \]
where $(\cdot,\cdot)$ denotes intersection product on $\mathrm{Div}(X)$. Our
task is thus to show that:
\[ 4g(g-1) (W_i-W_j,W_k) + (2g-1)(W_i-W_j,V_k) + (V_i-V_j,W_k)
= 2g(g-1) \nu(\ell_{ijk}) \, . \]
Let $S ' \to S$ be any finite cover of $S$, and let $X' \to S'$ be the minimal
desingularisation of the base change of $X \to S$ along $S' \to S$. By
functoriality of $(Q_{ik},q_{ik})$ and invariance under pullback 
it suffices to prove the formula for $X' \to
S'$. We start with the case that the residue characteristic of $R$ is not equal
to~$2$. Let $m$ be the maximal ideal of $R$ and let $K$ be the fraction field of
$R$. By \cite[Lemma 4.1]{ka} we may assume that $X \otimes K$ has an affine
equation $y^2 = \prod_{r=1}^{2g+2} (x-a_r)$ such that:
\begin{itemize}
\item the $a_r$ are distinct elements of $R$;
\item the valuations $\nu(a_r-a_s)$ are even for $r \neq s$;
\item the $a_r$ lie in at least~$3$ distinct residue classes of $R$ modulo $m$.
\end{itemize}
As the sections $W_i,W_j,W_k$ are disjoint we are reduced to showing that:
\[ (2g-1)(W_i-W_j,V_k) + (V_i-V_j,W_k) = 2g(g-1) \nu(\ell_{ijk}) \, . \]
Let $\alpha$ be the subset $\{ a_1,\ldots,a_{2g+2} \}$ of $R$. To $\alpha$ we
associate a finite tree $T$, as follows: for each positive integer $n$ let
$\rho_n \colon \alpha \to R/m^n$ be the canonical residue map. Define
$\Lambda_n$ to be the set of residue classes $\lambda$ in $R/m^n$ such that
$\rho_{n}^{-1}(\lambda) \subset \alpha$ has at least~$2$ elements. The vertices
of $T$ are then the elements $\lambda$ of $\Lambda_n$ for $n$ running through the
non-negative integers; there are only finitely many such $\lambda$. The edges of
$T$ are the pairs $(\lambda,\lambda')$ of vertices $\lambda,\lambda'$ where
$\lambda \in \Lambda_n$, $\lambda' \in \Lambda_{n+1}$ and $\lambda' \mapsto
\lambda$ under the natural map $\Lambda_{n+1} \to \Lambda_n$, for some $n$. If
$\lambda$ is a vertex of $T$ there is a unique $n$ such that $\lambda \in
\Lambda_n$; we call $n$ the \emph{level} of $\lambda$. 

Let $F$ be the special
fiber of $\pi$ and let $\Gamma$ be the dual graph of $F$. 
According to \cite[Section~5]{bo} or \cite[Section~4]{ka} 
there is a natural graph morphism $\varphi \colon \Gamma \to
T$. If $C$ is a vertex of $\Gamma$ we denote by $\lambda_C$ the image  
of $C$ in
$T$, and by $n_C$ the level of $\lambda_C$. For each $r=1,\ldots,2g+2$ we denote
by $C_r$ the unique irreducible component of $F$ through which the
$\sigma$-invariant section of $\pi$ corresponding to $a_r$ passes, by
$\lambda_r$ the image of $C_r$ in $T$, and by $n_r$ the level of $\lambda_r$. 
By construction of $\varphi$ the element 
$a_r$ is a representative of the class $\lambda_r$ and
$ n_r = \max_{s \neq r} \nu(a_r - a_s) $. 

For each irreducible component $C$ of $F$ we choose a representative $a_C$
of $\lambda_C$ in $\alpha$. If $f$ is a non-zero rational function on $X$ we
denote by $\nu_C(f)$ the multiplicity of $f$ along $C$. We have:
\[ \nu_C(x - a_r) = \min \{ n_C, \nu(a_C-a_r)\} \, , \]
independent of the choice of $a_C$, by \cite{ka}, proof of Lemma~$5.1$. In
particular:
\[ \nu_{C_k}(x-a_r) = \left\{ \begin{array}{ll}
\nu(a_k - a_r) & r \neq k \\ n_k & r=k \end{array} \right. \]
and so:
\[ \nu_{C_k}(y) = \frac{1}{2} n_k + \frac{1}{2} \sum_{r \neq k} \nu(a_k -a _r)
\, . \]
Now write $V_k = \sum_C \mu_k(C) \cdot C$ with $\mu_k(C) \in \mathbb{Z}$ and
with $C$ running through the irreducible components of $F$. 
We claim that:
\[ \mu_k(C) = (g-1) \min \{ n_C, \nu(a_k-a_C) \} - \nu_C(y) + n_C -
(g-\frac{1}{2})n_k + \frac{1}{2} \sum_{r \neq k} \nu(a_k - a_r)   \]
for all $C$. To see this, consider the rational section:
\[ \omega_k = (x-a_k)^{g-1} \frac{\mathrm{d}x}{y} \]
of $\omega$. Let $\mu_k'(C) = (g-1)\nu_C(x-a_k) - \nu_C(y) + n_C$. According to
\cite[Lemma 5.2]{ka} we have:
\[ \mathrm{div}_X \, \omega_k = (2g-2)W_k + \sum_C \mu_k'(C) \cdot C \, . \]
It follows that $V_k$ and $\sum_C \mu_k'(C) \cdot C$ differ by a multiple of
$F$. As $W_k^* V_k$ is trivial we find that:
\[ V_k = \sum_C ( \mu_k'(C) - \mu_k'(C_k)) \cdot C \, . \]
As we have:
\begin{align*} \mu_k'(C_k) &=  (g-1) \nu_{C_k} (x-a_k) - \nu_{C_k}(y) + n_k
\\
&=  (g-1)n_k - \frac{1}{2}n_k - \frac{1}{2} \sum_{r \neq k} \nu(a_r-a_k) + n_k
\\
&= (g - \frac{1}{2})n_k - \frac{1}{2} \sum_{r \neq k} \nu(a_r -a_k) \, , 
\end{align*}
the claim follows. As an immediate consequence we have:
\begin{align*}
(W_i,V_k) & =  \mu_k(C_i) \\
& =  (g-1) \min \{ n_i, \nu(a_k-a_i) \} - \nu_{C_i}(y) + n_i -
(g-\frac{1}{2})n_k + \frac{1}{2} \sum_{r \neq k} \nu(a_k-a_r) \\
& = (g-1) \nu(a_k-a_i) + \frac{1}{2}n_i - \frac{1}{2} \sum_{r \neq i} \nu(a_r
-a_i) - (g-\frac{1}{2})n_k + \frac{1}{2} \sum_{r \neq k} \nu(a_k-a_r) \, , 
\end{align*}
so that:
\[ (W_i-W_j,V_k) = (g-1) \nu \left( \frac{a_i -a_k}{a_j-a_k} \right) +
\frac{1}{2} (n_i-n_j) + \frac{1}{2} \sum_{r \neq i,j} \nu \left(
\frac{a_j-a_r}{a_i-a_r} \right) \, , \]
and in a similar fashion:
\[ (V_i-V_j,W_k) = (g-1) \nu \left( \frac{a_i-a_k}{a_j-a_k} \right) - (g -
\frac{1}{2})(n_i-n_j) - \frac{1}{2} \sum_{r \neq i,j} \nu \left(
\frac{a_j-a_r}{a_i-a_r} \right) \, . \]
This leads to:
\[ (2g-1)(W_i-W_j,V_k) + (V_i-V_j,W_k) = 2g(g-1) 
\nu \left( \frac{a_i-a_k}{a_j-a_k} \right) + (g-1) \sum_{r \neq i,j} 
\nu \left( \frac{a_j-a_r}{a_i-a_r} \right) \, , \]
and the required formula follows.

Next we consider the case that $R$ does have residue characteristic equal 
to~$2$. We have that $X \to S$ is smooth, and we may assume 
that all Weierstrass points
of $X \otimes K$ are rational over $K$, hence extend to sections of $\pi$. The
divisors $V_i,V_j$ and $V_k$ are empty, and we are reduced to showing that
simply:
\[ 2(W_i-W_j,W_k) = \nu(\ell_{ijk}) \, . \]
Let $h \colon X \to Y$ be the quotient of $X$ by $\sigma$. According to
\cite[Section 5]{lk} we have that $Y \to S$ is a smooth proper family of curves
of genus~$0$, and $h$ is finite flat of degree~$2$. Let $P_i,P_j,P_k \colon S
\to Y$ denote the sections $W_i,W_j,W_k$, composed with $h$. By the projection
formula we have:
\[ 2(W_i,W_k) = (P_i,P_k) \, , \quad  2(W_j,W_k) = (P_j,P_k) \, . \]
By \cite[Lemma 6.1]{ka} we may assume that on an affine open subset $X$ is given
by an equation $y^2 + p(x)y = q(x)$ with $p,q \in R[x]$ such that
$p^2 + 4q$ is a separable polynomial of degree $d=2g+2$. As $p=0$ defines the
fixed point subscheme of the hyperelliptic involution on the special fiber of $X
\to S$ the coefficients of $p$ generate the unit ideal in $R$. It follows that
we may even assume after a translation that $p(0)$ is a unit in $R$ and subsequently
after making a coordinate transformation $ x \mapsto x^{-1}$ that $p$ has
degree~$g+1$ and leading coefficient a unit in $R$. Write $f=p^2 +4q \in R[x]$.
Then $f$ has leading coefficient a unit in $R$ as well. We can write
$f = b\cdot \prod_{r=1}^{2g+2} (x-a_r)$ in $K[x]$; then $b \in R^\times$, and 
by Gauss's Lemma the $a_i$ are actually in $R$. Let
$a_i,a_j,a_k \in R$ correspond to $P_i,P_j,P_k$. As $y^2 = f(x)$ is a
hyperelliptic equation for $X \otimes K$ we have:
\[ \ell_{ijk} = \frac{a_i-a_k}{a_j-a_k} \sqrt[2g]{ \prod_{r \neq i,j} 
\frac{a_j-a _r }{a_i -a_r }} = 
\frac{a_i-a_k}{a_j-a_k} \sqrt[2g]{-\frac{ f'(a_j)}{f'(a_i)} } \]
in $K^\times/\mu_{2g}$. Since by the projection formula: 
\[ 2(W_i,W_k) = \nu(a_i-a_k) \, , \quad 2(W_j,W_k) = \nu(a_j-a_k) \, , \]
we are done once we prove that $f'(a_j)/f'(a_i)$ is a unit in $R$. Let $a_r$ be
an arbitrary root of $f$; we will show that $\nu(f'(a_r))=\nu(4)$, so that
$\nu(f'(a_r))$ is
independent of $r$. From the equation $p(a_r)^2 + 4q(a_r) =0$ we obtain
first of all that $p(a_r)$ is divisible by~$2$ in $R$. 
From $f=p^2+4q$ in $R[x]$ we obtain 
$f'(a_r) = 2p(a_r)p'(a_r) + 4q'(a_r) $ so that~$4$ divides $f'(a_r)$ in $R$ and
hence $\nu(f'(a_r)) \geq
\nu(4)$. According to \cite[Proposition 6.3]{ka} we have however:
\[ \sum_{s=1}^{2g+2} \nu(f'(a_s)) = \nu(b^d \cdot \prod_{s \neq t} (a_s-a_t))
= \nu(b^{2d-2} \cdot \prod_{s \neq t} (a_s-a_t)) = (2g+2)\nu(4) \, , \]
and we conclude that $\nu(f'(a_r))=\nu(4)$ as required.
\end{proof}
Let $\threep$ be the moduli stack, over $\mathbb{Z}[1/2]$, of pairs $(X
\to S, (W_i,W_j,W_k))$ with $S$ a $\mathbb{Z}[1/2]$-scheme, 
$X \to S$ a smooth hyperelliptic curve of genus~$g$, 
and $(W_i,W_j,W_k)$ a triple of distinct $\sigma$-invariant sections of $X
\to S$. 
\begin{lem} The moduli stack $\mathcal{I}^*_{g,3}$ is irreducible, and smooth 
over $\mathrm{Spec}(\mathbb{Z}[1/2])$.
\end{lem}
\begin{proof} Let $\mathcal{I}_g$ be the moduli stack of smooth hyperelliptic
curves of genus $g$ over $\mathbb{Z}[1/2]$, 
and let $\mathcal{U}_1$ be the universal hyperelliptic curve over
$\mathcal{I}_g$. Define inductively for $n \geq 2$ 
the algebraic stack $\mathcal{U}_n$ over $\mathcal{U}_{n-1}$ as the base change of
$\mathcal{U}_1 \to \mathcal{I}_g$ along $\mathcal{U}_{n-1} \to \mathcal{I}_g$. In
particular, the algebraic stack $\mathcal{U}_n$ is a smooth hyperelliptic curve
over $\mathcal{U}_{n-1}$, for all $n \geq 2$. Let $\mathcal{V}_n$ be the fixed
point substack of the hyperelliptic involution of $\mathcal{U}_n$ over
$\mathcal{U}_{n-1}$; then $\mathcal{V}_n \to \mathcal{U}_n$ is a closed immersion,
and the induced map $\mathcal{V}_n \to \mathcal{U}_{n-1}$ is finite \'etale by
\cite[Corollary 6.8]{lk}. We are interested in the cases $n=1,2,3$. Put:
\[ \mathcal{A} = \mathcal{V}_2 \times_{\mathcal{U}_1} \mathcal{V}_1 \, , \quad
\mathcal{B} = \mathcal{V}_3 \times_{\mathcal{U}_2} \mathcal{V}_2 \, , \quad
\mathcal{C} = \mathcal{A} \times_{\mathcal{V}_2} \mathcal{B} \, . \]
As $\mathcal{U}_1$ is naturally the moduli stack of $1$-pointed hyperelliptic
curves over $\mathbb{Z}[1/2]$, we have a natural map $\threep \to \mathcal{U}_3$.
It factors via the closed immersion $\mathcal{C} \to \mathcal{U}_3$, and the
induced map $\threep \to \mathcal{C}$ is an open immersion. Since $\mathcal{C} \to
\mathcal{I}_g$ is \'etale and the structure map $\mathcal{I}_g \to
\mathrm{Spec}(\mathbb{Z}[1/2])$ is smooth (see \cite[Theorem 3]{ll}), 
we obtain that $\threep$ is smooth
over $\mathrm{Spec}(\mathbb{Z}[1/2])$. It follows that
the generic points of $\threep$ all lie above the
generic point of $\mathrm{Spec}(\mathbb{Z}[1/2])$. But $\threep \otimes
\mathbb{Q}$ is a quotient of the moduli stack of $2g+2$ distinct points on
$\mathbb{P}^1$, which is irreducible, hence $\threep \otimes \mathbb{Q}$ is
irreducible. It follows that $\threep$ is irreducible as well.
\end{proof}
We need a suitable compactification of $\threep$. Note that we have a
natural closed immersion $\threep \to \mathcal{M}_{g,3}$ where $\mathcal{M}_{g,3}$
is the moduli stack, over $\mathbb{Z}[1/2]$, of $3$-pointed smooth projective
curves of genus~$g$. We let $\clos$ be the stack-theoretic closure of $\threep$
in $\overline{\mathcal{M}}_{g,3}$, the Knudsen-Mumford compactification of
$\mathcal{M}_{g,3}$. Then $\clos$ is integral, and 
proper and flat over
$\mathbb{Z}[1/2]$. We have $\ell_{ijk}^{2g}$ and, via
$\psi_{ijk}$, also $q_{ik}^{-1} \otimes q_{jk}$ as rational sections of the
structure sheaf of $\clos$. In fact, for every generically smooth semistable
hyperelliptic curve $\pi \colon X \to S$ with $S$ an integral
$\mathbb{Z}[1/2]$-scheme together with a triple of distinct $\sigma$-invariant 
sections with image in
the smooth locus of $\pi$, we have a natural period map $S \to \clos$ such that
both $q_{ik}^{-1} \otimes q_{jk}$ and $\ell_{ijk}^{2g}$ associated to $\pi$ on
$S$ are obtained by
pullback from $\clos$. 

It suffices therefore to prove Theorem \ref{main} 
for the tautological curve over $\clos$.
\begin{proof}[Proof of Theorem \ref{main}]
We first prove that $q_{ik}^{-1} \otimes q_{jk}$ and $\ell_{ijk}^{2g(g-1)}$
on $\clos$ differ by a unit in $\mathbb{Z}[1/2]$. Let $\nu \colon \norm \to \clos$ be the
normalisation of $\clos$. As $\mathbb{Z}[1/2]$ is an excellent ring, the
morphism $\nu$ is finite birational. In particular $\norm$ and $\clos$ are
isomorphic over an open dense substack. Hence, in order to prove that 
$q_{ik}^{-1} \otimes q_{jk}$ and $\ell_{ijk}^{2g(g-1)}$
differ by a unit in $\mathbb{Z}[1/2]$ it suffices to prove that their pullbacks
along $\nu$ do so over $\norm$. This goes in two steps: first we prove that 
$q_{ik}^{-1} \otimes q_{jk}$ and $\ell_{ijk}^{2g(g-1)}$
differ by an invertible regular function on $\norm$, and then that the set of
such functions is precisely $\mathbb{Z}[1/2]^\times$. As to the first step, since
$\norm$ is normal, it suffices to prove that for every point $x'$ of $\norm$ of
height~$1$, the sections $q_{ik}^{-1} \otimes q_{jk}$ and $\ell_{ijk}^{2g(g-1)}$
differ by a unit in the local ring $\mathcal{O}_{x'}$ at $x'$. So let
$x'$ be a point of height~$1$ on $\norm$. Let $\mathrm{Spec} (\mathcal{O}_{x'})
\to \norm \xrightarrow{\nu} \clos$ be the canonical map. Then the generic point
of $\mathrm{Spec} (\mathcal{O}_{x'})$ maps to the generic point of $\clos$ and by
pullback we obtain a generically smooth stable hyperelliptic curve $\pi' \colon
X' \to \mathrm{Spec} (\mathcal{O}_{x'})$ over $\mathcal{O}_{x'}$ 
together with a triple of distinct
$\sigma$-invariant sections with image in the smooth locus of $\pi'$. By
functoriality of $(Q_{ik},q_{ik})$ in passing from $X'$ to a minimal desingularisation of $X'$ over $\mathrm{Spec}
(\mathcal{O}_{x'})$ we may assume for our purposes that $X'$ is itself regular and
semistable over $\mathrm{Spec} (\mathcal{O}_{x'})$. Proposition \ref{firststep}
then gives us precisely what we need. 

Next let $p \colon \norm \to \mathrm{Spec} (\mathbb{Z}[1/2])$ be the structure
map; we claim that $p_* \mathcal{O}_{\norm}$ is equal to 
$ \mathcal{O}_{\mathrm{Spec}
(\mathbb{Z}[1/2])}$. This will prove that the set of invertible regular functions
on $\norm$ is equal to $\mathbb{Z}[1/2]^\times$. The claim follows from: (i) $p$ is
proper, and (ii) $p$ is flat and $\norm \otimes \mathbb{Q}$ is irreducible.
Indeed, if $p$ is proper and flat then $p_* \mathcal{O}_{\norm}$ is a finite
torsion-free $\mathcal{O}_{\mathrm{Spec} (\mathbb{Z}[1/2])}$-module. If moreover
$\norm \otimes \mathbb{Q}$ is irreducible then $p_* \mathcal{O}_{\norm}$ has generic
rank~$1$, and $p_* \mathcal{O}_{\norm} = \mathcal{O}_{\mathrm{Spec}
(\mathbb{Z}[1/2])}$ follows. That (i) is satisfied is clear as both 
$\nu$ and $\clos
\to \mathrm{Spec} (\mathbb{Z}[1/2])$ are proper. To see (ii), recall that
$\clos$ is irreducible and has its generic point mapping to the generic point
of $\mathrm{Spec} (\mathbb{Z}[1/2])$. The same holds for $\norm$ since $\nu$
is birational; in particular $\norm$ is flat over~$\mathbb{Z}[1/2]$.

The next step is to prove that the unit $u$ in $\mathbb{Z}[1/2]$ connecting
$q_{ik}^{-1} \otimes q_{jk}$ and $\ell_{ijk}^{2g(g-1)}$ is either $+1$ or $-1$.
For this we use a smooth hyperelliptic curve $X \to S$ 
with a triple of $\sigma$-invariant sections over the spectrum $S$ of a discrete
valuation ring with residue characteristic equal to~$2$, and generic
characteristic zero (such data exist). We have a period map $S
\otimes \mathbb{Q} \to
\clos \otimes \mathbb{Q}$ and applying once more Proposition \ref{firststep} we see that the
exponent of~$2$ in $u$ is vanishing. 

We finish by proving that $u=(-1)^{g-1}$. Cyclic permutation of the three 
$\sigma$-invariant sections in the moduli data induces a group of automorphisms 
of the moduli stack 
$\mathcal{I}^*_{g,3}$ of
order three. Its action on the regular functions on $\mathcal{I}^*_{g,3}$ 
yields the regular functions $\ell_{jki}^{2g}$ and $\ell_{kij}^{2g}$ 
from $\ell_{ijk}^{2g}$, as well as the functions
$q_{ji}^{-1} \otimes q_{ki}$ and $q_{kj}^{-1} \otimes q_{ij}$ 
from $q_{ik}^{-1} \otimes q_{jk}$. A small computation shows 
the cocyle relation:
\[ \ell_{ijk}^{2g} \cdot \ell_{jki}^{2g} \cdot \ell_{kij}^{2g} = -1 \, , \]
whereas we have:
\[ q_{ik}^{-1} \otimes q_{jk} \otimes q_{ji}^{-1} \otimes q_{ki} \otimes 
q_{kj}^{-1} \otimes q_{ij} = 1  \]
by the canonical symmetry isomorphism
(\ref{symmetry}). Now write
$\ell_{ijk}^{2g(g-1)} = u \cdot q_{ik}^{-1} \otimes q_{jk}$. As $u$, being a constant
function, is left invariant by any automorphism of $\mathcal{I}^*_{g,3}$, 
we obtain the identities:
\[ \ell_{jki}^{2g(g-1)} = u \cdot q_{ji}^{-1} \otimes q_{ki} \, , \quad
\ell_{kij}^{2g(g-1)} = u \cdot q_{kj}^{-1} \otimes q_{ij} \, . \]
Combining the cocycle relations with these identities 
yields $u=u^3 = (-1)^{g-1}$.
The proof of Theorem \ref{main} is thereby complete.
\end{proof}
We are now ready to prove Theorem C. Let $K$ be a field
which is either a complete discrete valuation field, or $\mathbb{R}$ or
$\mathbb{C}$. Let $\bar{K}$ be an algebraic closure of $K$ and endow $\bar{K}$
with its canonical absolute value $|\cdot|$. This absolute value is defined as
follows: if $K$ is
a complete discrete valuation field, endow $K$ with the absolute value
$|\cdot|_K$
such that $|\pi|=1/e$ for a uniformizer $\pi$ of $K$; we get an absolute value
$|\cdot|$ on $\bar{K}$ by taking the unique extension of $|\cdot|_K$ to $\bar{K}$. If
$K=\mathbb{R}$ or $K =\mathbb{C}$ we endow $\bar{K}=\mathbb{C}$ with the
standard euclidean norm. 

Assume that $K$ does not have characteristic~$2$.
\begin{thm} \label{thmC} 
Let $X$ be a hyperelliptic curve over $K$, and let $(,)_a$ be Zhang's 
admissible pairing on $\mathrm{Div}(X \otimes \bar{K})$. Let $w_i,w_j,w_k$ be three distinct Weierstrass points
on $X \otimes \bar{K}$. Then the formula:
\[ (w_i-w_j,w_k)_a = -\frac{1}{2} \log | \ell_{ijk} | \]
holds.
\end{thm}
\begin{proof} We apply Theorem \ref{main} to the morphism $X \otimes
\bar{K} \to \mathrm{Spec} (\bar{K})$. Let $g$ be the genus of~$X$. Under the
isomorphism $\psi_{ijk}$ the sections $q_{ik}^{-1} \otimes q_{jk}$ and
$(-\ell_{ijk}^{2g})^{g-1}$ are identified, and by Proposition \ref{Qik} the
isomorphism $\psi_{ijk}$ is an admissible isometry. Let $|\cdot|_a$ be the
admissible norm on $-Q_{ik} +Q_{jk}$ on $\mathrm{Spec} (\bar{K})$. We obtain:
\begin{align*}
2g(g-1) \log | \ell_{ijk} | &= \log | q_{ik}^{-1} \otimes q_{jk} |_a \\
 & =  -4g(g-1) (w_i-w_j,w_k)_a  \, , 
\end{align*}
and the theorem follows.
\end{proof}
\begin{cor} Let $w_i,w_j,w_k,w_r$ be four distinct Weierstrass points on $X
\otimes \bar{K}$. Let $\mu_{ijkr}$ in $\bar{K}$ be the cross-ratio on
$w_i,w_j,w_k,w_r$. Then the formula:
\[ (w_i-w_j,w_k-w_r)_a = -\frac{1}{2} \log |\mu_{ijkr}| \]
holds.
\end{cor}
\begin{proof} This follows directly from the identity:
\[ \frac{ \ell_{ijk} }{\ell_{ijr} } = \frac{ a_i - a_k }{a_j - a_k} \cdot \frac{
a_j-a_r}{a_i - a_r} = \mu_{ijkr} \]
which is easily checked.
\end{proof}
\begin{remark} In the case that $K =\mathbb{R}$ or $\mathbb{C}$ the admissible
pairing is given by the Arakelov-Green's function $G$ of the compact Riemann surface
$X(\bar{K})$. Theorem \ref{thmC} translates into the remarkable formula:
\[ \frac{G(w_i,w_k)}{G(w_j,w_k)} = \sqrt{|\ell_{ijk}|} \, . \]
It would be interesting to see if one could give a direct proof of this formula
that does not use moduli spaces.
\end{remark}

\section{The invariant $\chi$} \label{chi}

In this section we introduce the $\chi$-invariant as announced 
in the Introduction. The definition may seem rather \emph{ad hoc} at first sight, but
in the function field context the invariant already occurs, as mentioned before,
in work of A. Moriwaki \cite{mo} and K. Yamaki \cite{ya2}.
In the next section we
present a more intrinsic approach to $\chi$, using the arithmetic of 
symmetric roots.

Let $\pi \colon X \to S$ be an
arbitrary generically smooth semistable hyperelliptic curve of genus $g$ with
$S$ an integral scheme. Let
$\omega$ be the relative dualising sheaf of $\pi$, and let $\lambda = \det
R\pi_* \omega$ be the Hodge bundle on $S$. As is 
explained in \cite[Section 2]{ma}, the line bundle $(8g+4)\lambda$ has a
canonical non-zero rational section $\Lambda_g$, satisfying the following
properties. We write $\mathcal{I}_g$ for the moduli stack of smooth hyperelliptic
curves of genus $g$ over $\mathbb{Z}$. 
\begin{itemize}
\item the formation of $\Lambda_g$ is compatible with dominant base change;
\item if $S$ is normal, then $\Lambda_g$ is global;
\item if $\pi$ is smooth, then $\Lambda_g$ is global and nowhere vanishing; 
\item if $S$ is
the spectrum of a field $K$ which is not of characteristic~$2$, and
$y^2+p(x)y=q(x)$ is an affine equation of $X$ with $p,q \in K[x]$, one has:
\[ \Lambda_g = (2^{-(4g+4)} \cdot D)^g \cdot \left( \frac{ \mathrm{d} x}{2y+p}
\wedge \cdots \wedge \frac{ x^{g-1} \mathrm{d} x}{2y + p} \right)^{\otimes 8g+4}
\, , \]
where $D$ is the discriminant of the separable polynomial $p^2+4q \in K[x]$; 
\item let $\overline{\mathcal{I}}_g$ be the stack-theoretic closure of $\mathcal{I}_g$
in $\overline{\mathcal{M}}_g$, the moduli stack of stable curves of genus $g$
over $\mathbb{Z}$. Then on $S$, the rational section $\Lambda_g$ can be 
obtained as the pullback of a rational section $\Lambda_g$ of $(8g+4)\lambda$ 
on $\overline{\mathcal{I}}_g$. 
\end{itemize}

Now assume that $S = \mathrm{Spec}(K)$ where $K$ is either a complete discrete
valuation field or $\mathbb{R}$ or $\mathbb{C}$. Note that we put no
restrictions, at this stage, on the characteristic of $K$. The section 
$\Lambda_g$ gives rise to a real-valued invariant $d(X)$ associated to $X$, as
follows. If $K$ is non-archimedean, let $\mathcal{X} \to \mathrm{Spec}(R)$ be
the minimal regular model of $X$, where $R$ is the valuation ring of $K$. Then
$d(X)$ is the order of vanishing of $\Lambda_g$ along the closed point of $S$.
If $K$ is archimedean, the $\mathbb{C}$-vector space 
$\mathrm{H}^0(X(\mathbb{C}),\omega)$ is equipped with a natural hermitian inner 
product $(\omega,\eta) \mapsto \frac{i}{2} \int_{X(\mathbb{C})} \omega \, \bar{\eta}$,
and $d(X)$ is the $-\log$ of the norm of $\Lambda_g$ with
respect to this inner product. 

One can give explicit formulas for $d(X)$. Start again 
with the case that $K$ is a
complete discrete valuation field. Assume that $X$ has semistable reduction 
over~$K$, and let as above $\mathcal{X} \to \mathrm{Spec}(R)$ be the minimal regular
model of $X$. Let $x$ be a singular point in the special fiber of $\mathcal{X}$.
We say that $x$ is of type $0$ if the local normalisation of the special fiber
at $x$ is connected, and that $x$ is of type $i$, where 
$1 \leq i \leq \lfloor g/2 \rfloor$, if the local 
normalisation of the special fiber of $\mathcal{X}$ at $x$ is the disjoint 
union of two semistable curves of genus $g$ and $g-i$.  
Let $x$ be a singular point of type $0$. Let $\sigma$ be the hyperelliptic 
involution of $\mathcal{X}$. We have the following two
possibilities for $x$:
\begin{itemize}
\item $x$ is fixed by $\sigma$. Then we say $x$ is of subtype $0$;
\item $x$ is not fixed by $\sigma$. Then the local normalisation of the special
fiber of $\mathcal{X}$ at $\{ x,\sigma(x) \}$ consists of two connected
components of genus $j$ and $g-j-1$, say, where $1 \leq j \leq \lfloor (g-1)/2
\rfloor$. In this case we say that $x$ is of subtype $j$. 
\end{itemize}
Let $\delta_i(X)$ for $i=1,\ldots,\lfloor g/2 \rfloor$ be the number of
singular points in the special fiber of $\mathcal{X}$ 
of type $i$, let $\xi_0(X)$ be the number
of singular points of subtype $0$, and let $\xi_j(X)$ for $j=1,\ldots,\lfloor
(g-1)/2 \rfloor$ be the number of pairs of nodes of subtype $j$. Then the
following equality holds, proved in increasing order of generality by
M.~Cornalba and J.~Harris, I.~Kausz, and K.~Yamaki \cite{ya1}:
\begin{equation} \label{ch} 
d(X) = g \xi_0(X) + \sum_{j=1}^{ \lfloor (g-1)/2 \rfloor} 2(j+1)(g-j)
\xi_j(X) + \sum_{i=1}^{\lfloor g/2 \rfloor} 4i(g-i) \delta_i(X) \, . 
\end{equation}
If $K$ equals $\mathbb{R}$ or $\mathbb{C}$ then $d(X)$ can be related to a product of
Thetanullwerte, as explained in \cite[Section 8]{dj1}. Let $\tau$ in the Siegel
upper half space be a normalised period matrix for $X(\mathbb{C})$ 
formed on a canonical symplectic
basis of $\mathrm{H}_1(X(\mathbb{C}),\mathbb{Z})$. Let $\varphi_g$ be the
level $2$ Siegel modular form from \cite[Definition 3.1]{lo} and put:
\[ \Delta_g = 2^{-(4g+4)n} \varphi_g \, , \]
where $n = {2g \choose g+1}$. Then the real number $\| \Delta_g \|(X) = (\det
\mathrm{Im} \tau)^{2r} |\Delta_g(\tau)|$ is independent of the choice of $\tau$,
where $r = {2g+1 \choose g+1}$, hence defines an invariant of $X(\mathbb{C})$.
It follows from \cite[Proposition 3.2]{lo} that the formula:
\begin{equation} \label{archimd} n  d(X) = -n \log \| \Lambda_g \| = -4g^2r\log(2\pi) - g \log 
\|\Delta_g \|(X) 
\end{equation}
holds. 

The invariant $\chi(X)$ of $X$ is determined by $d(X)$ and the invariants
$\varepsilon(X)$ and $\delta(X)$ which we discuss next.  The invariant
$\varepsilon(X)$ stems from \cite{zh1} and is defined as follows.  Let $K$ again
be non-archimedean. We keep the assumption that $X$ has semistable reduction
over $K$.  Let $R(X)$ be the reduction graph of $X$, let $\mu_X$ be the
admissible measure on $R(X)$, let $K_X$ be the canonical divisor on $R(X)$, and
let $g_X$ be the admissible Green's function on $R(X)$. Then: 
\[ \varepsilon(X)
= \int_{R(X)} g_X(x,x)((2g-2)\mu_X + \delta_{K_X}) \, . \] 
If $K$ is archimedean, one simply puts $\varepsilon(X)=0$. 

As to $\delta(X)$, for $K$ non-archimedean we put
$\delta(X) = \sum_{i=0}^{\lfloor g/2 \rfloor} \delta_i(X)$, the number of
singular points on the special fiber of the minimal regular model of $X$ over
$K$; if $K$ equals $\mathbb{R}$ or $\mathbb{C}$ we put $\delta(X)=-4g
\log(2\pi) + \delta_F(X)$ where $\delta_F(X)$ is the Faltings delta-invariant
of $X(\mathbb{C})$ defined on \cite[p.~402]{fa}.

The invariant $\chi(X)$ is determined by the following equality:
\begin{equation} \label{defchi} 
(2g-2) \chi(X) = 3 d(X) - (2g+1)(\varepsilon(X) + \delta(X)) \, . 
\end{equation}
Let $K$ be non-archimedean. It is clear that $\chi(X)=0$ if $X$ has good
reduction, since each of $d,\varepsilon$ and $\delta$ vanishes in this case.

From \cite{mo} one can calculate $\chi(X)$ in the case that $g=2$, based on the
classification of the semistable fiber types in genus~$2$. 
We display the results in a table:
\begin{center}
\begin{tabular}{|l|cccc|}
\hline  
Type & $d/2$ & $\delta $ & $\varepsilon$ & $\chi$ \\  [2pt] 
\hline  
$I$  & $0$ & $0$ & $0$ & $0$ \\  [4pt]
$II(a)$ & $2a$ & $a$ & $a$ & $a$ \\  [4pt] 
$III(a)$ & $a$ & $a$ & $\frac{1}{6}a$ & $\frac{1}{12}a $ \\ [4pt]
$IV(a,b)$ & $2a+b$ & $a+b$ & $a+\frac{1}{6}b$ & $a + \frac{1}{12}b$ \\ [4pt]
$V(a,b)$ & $a+b$ & $a+b$ & $\frac{1}{6}(a+b)$ & $ \frac{1}{12}(a+b)$ \\ [4pt]
$VI(a,b,c)$ & $2a+b+c$ & $a+b+c$ & $a+ \frac{1}{6}(b+c)$ & $a + \frac{1}{12}(b+c) $ \\
[4pt]
$VII(a,b,c)$ & $a+b+c$ & $a+b+c$ & 
$\frac{1}{6}(a+b+c) + \frac{1}{6} \frac{abc}{ab+bc+ca}$ &
$\frac{1}{12}(a+b+c) - \frac{5}{12}\frac{abc}{ab+bc+ca}$ \\  [4pt]
\hline
\end{tabular}
\end{center}
\vspace{8pt} 
In the case $g \geq 3$ one has an effective lower bound for $\chi(X)$
which is strictly positive in the case of non-smooth reduction, by work of
Yamaki \cite{ya2}. We quote his result:
\[ \chi(X) \geq \frac{ (2g-5)}{24g} \xi_0(X) + \sum_{j=1}^{\lfloor (g-1)/2
\rfloor} \frac{3j(g-1-j)-g-2 }{ 3g} \xi_j(X) + \sum_{i=1}^{\lfloor g/2 \rfloor}
\frac{ 2i(g-i)}{g} \delta_i(X)  \]
if $g \geq 5$, and:
\[ \chi(X) \geq \frac{ (2g-5)}{24g} \xi_0(X) + \sum_{j=1}^{\lfloor (g-1)/2
\rfloor} \frac{2j(g-1-j)-1 }{ 2g} \xi_j(X) + \sum_{i=1}^{\lfloor g/2 \rfloor}
\frac{ 2i(g-i)}{g} \delta_i(X) \]
for the cases $g=3,4$. The most difficult part of the proof lies in obtaining 
suitable upper bounds for the $\varepsilon$-invariant ranging over 
the reduction graphs of hyperelliptic curves of a fixed genus using 
combinatorial optimisation. 

If $K$ is archimedean, one easily gets an exact formula for $\chi(X)$ from
(\ref{archimd}) and (\ref{defchi}). We state the result for completeness: 
\[ \chi(X) = -\frac{ 8g(2g+1)}{2g-2} \log(2\pi) - \frac{3g}{(2g-2)n} \log
\|\Delta_g\|(X) - \frac{2g+1}{2g-2} \delta_F(X) \, . \]

If $K$ is a number field or the function field of a curve over a field $k_0$, and $X$
has semistable reduction over $K$,
one has the following
global formulas involving $d,\varepsilon,\delta$ and the places $v$ of $K$. Let
$Nv$ be the following local factors at $v$: in the number field case and for $v$
non-archimedean, put $ Nv =  \# \kappa(v)$ with $\kappa(v)$ the residue
field at $v$; for $v$ archimedean put $\log Nv=1$ if $v$ is real, and $\log
Nv=2$ if $v$ is complex. In the function field case, let $\log Nv$ be the degree
of $v$ over the base field $k_0$. Then first of all: 
\[ (8g+4) \deg \det R\pi_* \omega = \sum_v d(X_v) \log Nv \, , \]
where $\deg \det R\pi_* \omega$ is the geometric degree of $\det R\pi_* \omega$
in the function field case, and the non-normalised Faltings height of $X$ in the
number field case. This follows directly from the definition of $d$. 
Next one has:
\[ (\omega,\omega)_a = (\omega,\omega) - \sum_v \varepsilon(X_v) \log Nv \]
by \cite[Theorem 4.4]{zh2}, where $(\omega,\omega)_a$ and $(\omega,\omega)$ are
the admissible and usual self-intersections of the relative dualising sheaf of
$X$ over $K$, respectively. Finally:
\[ 12 \deg \det R\pi_* \omega = (\omega,\omega) + \sum_v \delta(X_v) \log Nv \]
which is the Noether formula for $X$ over $S$, cf. \cite[Th\'eor\`eme 2.5]{mb}
for the number field case.
 
We easily find the formula: 
\begin{equation} \label{formulachi}  
(\omega,\omega)_a = \frac{2g-2}{2g+1} \sum_v \chi(X_v) \log Nv    
\end{equation}
expressing $(\omega,\omega)_a$ in terms of the $\chi(X_v)$.  From the results of
Moriwaki and Yamaki mentioned above one gets an effective version of the
Bogomolov conjecture for $X$, if $K$ is a function field.

\section{Intrinsic approach to $\chi$} \label{intrinsic}

In this section we provide an alternative approach to $\chi$. We
construct a canonical non-zero rational section $q$ of the line 
bundle $(2g+1) \langle
\omega,\omega \rangle$ on the base of a generically smooth semistable
hyperelliptic curve of genus $g$, and show that  
$\chi$ is essentially the $-\log$ of the admissible norm of $q$.

A crucial ingredient of our approach is the arithmetic of the 
symmetric discriminants of a hyperelliptic curve. 
These are intimately related to the curve's symmetric roots, 
see \cite[Section 2]{gu} for a discussion. The definition is as follows. 
Let $\kappa$ be a field not of characteristic~$2$ 
and let $X$ be a hyperelliptic
curve of genus $g \geq 2$ over $\kappa$. Let $\bar{\kappa}$ be a
separable algebraic closure of $\kappa$, and let $w_1,\ldots,w_{2g+2}$ be the
Weierstrass points of $X \otimes \bar{\kappa}$. Let $(w_i,w_j)$ be a pair of
these. We have well-defined sets of symmetric roots $\{
\ell_{ijk}\}^{\zeta}_{k \neq i,j}$ 
in $\bar{\kappa}$ associated to $i,j$, parametrised by the elements $\zeta$ of
$\mu_{2g}$. The $\zeta$'s give rise to a set of \emph{symmetric equations}:
\[ \mathcal{C}_{ij}^\zeta = x \prod_{k \neq i,j} (x - \ell_{ijk}) \]
for $X$ over $\bar{\kappa}$, parametrised by $\zeta$. The discriminant:
\[ d_{ij} = \prod_{r,s \neq i,j \atop r \neq s} 
(\ell_{ijr}  - \ell_{ijs} ) \]
of a $\mathcal{C}_{ij}^\zeta$ is independent of the choice of $\zeta$ and is
called the \emph{symmetric discriminant} of the pair $(w_i,w_j)$. It is a
well-defined element of the field of definition inside $\bar{\kappa}$ of 
$(w_i,w_j)$. A small calculation shows that the formula:
\begin{equation} \label{rootsfromdisc} 
\frac{d_{ik}}{d_{jk}} = -\ell_{ijk}^{2g(2g+1)}  
\end{equation}
holds in $\bar{\kappa}$ for all $k \neq i,j$, allowing us to compute suitable powers of 
the symmetric roots of $X$ from the symmetric discriminants of $X$.

We start our construction of the rational section $q$ with the case of a smooth hyperelliptic curve
$\pi \colon X \to S$ of genus $g$ where $S$ is a scheme whose generic
points are not of characteristic~$2$. By \cite[Proposition 7.3]{lk} there exists a
faithfully flat morphism $S' \to S$ such that the smooth hyperelliptic curve 
$X \times_S S' \to S'$ has $2g+2$ sections $W_1,\ldots,W_{2g+2}$ invariant for
the hyperelliptic involution. Let $(W_i,W_k)$ be a pair of these. As we saw in
Section \ref{admpairing} there exists a line bundle $Q_{ik}$ on $S'$ associated to
$(W_i,W_k)$ together with a canonical non-zero rational section $q_{ik}$ of
$Q_{ik}$. We can and will view $q_{ik}$ as a rational section
of $\langle \omega,\omega \rangle$ on $S'$, by Proposition \ref{Qik}. 

Let $d_{ik}$ be the symmetric discriminant of $(W_i,W_k)$, 
viewed as a rational function on $S'$. We define:
\begin{equation} \label{formulaq} 
q = (2^{4g} d_{ik} )^{g-1} \cdot q_{ik}^{\otimes 2g+1} \, , 
\end{equation}
viewed as a rational section of $(2g+1) \langle \omega, \omega \rangle$ on $S'$.
\begin{lem} \label{qforsmooth}  The rational section $q$ of $(2g+1) \langle
\omega, \omega \rangle$ is independent of the choice of $(W_i,W_k)$, and
descends to a canonical rational section of $(2g+1) \langle \omega, \omega
\rangle$ on $S$.  
\end{lem}
\begin{proof}
To see this, first fix an index $k$ and consider the sections
$(2^{4g} d_{ik} )^{g-1} \cdot q_{ik}^{\otimes 2g+1}$ and 
$(2^{4g} d_{jk} )^{g-1} \cdot q_{jk}^{\otimes 2g+1}$ for $i,j \neq k$. According to 
equation (\ref{rootsfromdisc}) we have:
\[ \frac{d_{ik}}{d_{jk}} = -\ell_{ijk}^{2g(2g+1)} \, , \]
whereas by Theorem \ref{main} we have:
\[ q_{ik}^{-1} \otimes q_{jk} = (-\ell_{ijk}^{2g})^{g-1} \, . \]
It follows that the $(2^{4g} d_{ik} )^{g-1} \cdot q_{ik}^{\otimes 
2g+1}$ are mutually
equal, where $i$ runs over the indices different from $k$. 
By symmetry considerations 
we can vary $k$ as well and the independence of $q$ on the choice of
$(i,k)$ follows. By faithfully flat descent, see \cite[Expos\'e VIII, Th\'eor\`eme
1.1]{sga1}, we obtain that $q$ comes from the
base $S$. 
\end{proof}
Let again $\mathcal{I}_g$ be the moduli stack of smooth hyperelliptic curves
of genus $g$ over~$\mathbb{Z}$. 
By Lemma \ref{qforsmooth} we have $q$ as a canonical 
rational section of the line bundle $(2g+1) \langle \omega,\omega \rangle$ on
$\mathcal{I}_g$, and by pullback we obtain $q$ on the base of \emph{any} smooth
hyperelliptic curve. Even better, 
by extension we get $q$ as a rational section on the base
of any \emph{generically smooth} hyperelliptic curve. We isolate this result in
a theorem.
Let $S$ be an integral scheme, let
$\pi \colon X \to S$ be a generically smooth hyperelliptic curve of genus $g\geq
2$ over $S$, and let $\omega$ be the relative dualising sheaf of $\pi$.
\begin{thm} The line bundle 
$(2g+1) \langle \omega,\omega \rangle$ on $S$ 
has a canonical rational section~$q$. If $S$ does not have generic
characteristic equal to~$2$, then $q$ is given by 
equation (\ref{formulaq}). The formation of $q$ is compatible 
with dominant base change.
\end{thm}
The next result yields Theorems A and B as an immediate consequence.
\begin{thm} \label{chithm} 
Assume that $S$ is a normal integral scheme. 
\begin{itemize}
\item[(i)] the rational section $q$ of $(2g+1)\langle \omega,\omega \rangle$ 
is in fact a global section, with no
zeroes if $\pi$ is smooth. 
\end{itemize}
Let $K$ be either a complete discrete valuation field or let $K$ equal 
$\mathbb{R}$ or $\mathbb{C}$. Let $\bar{K}$ be an algebraic closure of $K$ and assume that
$S=\mathrm{Spec} (\bar{K})$. Let $|\cdot|_a$
be the admissible norm on $(2g+1) \langle \omega,\omega \rangle$. Then:
\begin{itemize}
\item[(ii)] the formula: 
\[ -\log |q|_a = (2g-2) \chi(X) \]
holds;
\item[(iii)] if $K$ does not have characteristic~$2$, and  
$w_1,\ldots,w_{2g+2}$ on $X \otimes \bar{K}$ are the Weierstrass points of $X$,
then the formula:
\[ - \log |q|_a = -4g(g-1) \left( \log|2| + \sum_{k \neq i} (w_i,w_k)_a \right)
\]
holds, for all $i=1,\ldots,2g+2$, where $(,)_a$ denotes Zhang's admissible
pairing on $\mathrm{Div}(X \otimes \bar{K})$.
\end{itemize}
\end{thm}
\begin{proof} To prove that $q$ has no zeroes if $\pi$ is smooth, 
note that it suffices to prove this in the case of
the tautological curve over $\mathcal{I}_g$. We have that $\mathcal{I}_g$ is
normal, as $\mathcal{I}_g \to \mathrm{Spec} (\mathbb{Z})$ is smooth by
\cite[Theorem 3]{ll}. Thus it is sufficient to prove the statement for the
case of a smooth hyperelliptic curve
$\pi \colon X \to S$ where $S$ is the spectrum of a discrete valuation ring $R$ 
of characteristic zero. We can and will assume that all Weierstrass points of 
the generic fiber of $\pi$ are rational, hence extend to sections
$W_1,\ldots,W_{2g+2}$ of $\pi$. Let
$\nu(\cdot)$ denote order of vanishing along the closed point of $S$, and fix an
index $i$. It is not hard to check the formula $\prod_{k \neq i}
d_{ik}=1$ for the symmetric discriminants. This gives us, directly from the
definition of $q$:
\[
\nu(q) = 4g(g-1) \nu(2) + \sum_{k \neq i} \nu(q_{ik}) \, .  \]
Hence:
\[ \nu(q) = 4g(g-1) \left( \nu(2) - \sum_{k \neq i} (W_i,W_k) \right) \, , \]
where $(,)$ denotes intersection product on $\mathrm{Div}(X)$. Note that the
$V$'s are empty.
If the residue characteristic of $R$ is not equal to~$2$ we immediately obtain
the vanishing of $\nu(q)$ since $W_i$ is disjoint from each $W_k$. 
So assume that the
residue characteristic of $R$ is equal to~$2$. It follows from the proof of
Proposition \ref{firststep} that:
\[ 2 \sum_{k \neq i} (W_i,W_k) = \nu(4) \, . \]
We see that $\nu(q)$ vanishes in this case as well. 
This proves the second half of~(i).

Now let as above $\lambda=\det R\pi_* \omega$ be the Hodge bundle on
$\overline{\mathcal{I}}_g$, the stack-theoretic closure of $\mathcal{I}_g$ in
$\overline{\mathcal{M}}_g$, and let $\delta$ be the line bundle associated to
the restriction of the boundary divisor of $\overline{\mathcal{M}}_g$ to
$\overline{\mathcal{I}}_g$. 
By \cite[Th\'eor\`eme 2.1]{mb} there exists an isomorphism: 
\[ \mu \colon 3(8g+4) \lambda - (2g+1) \delta \xrightarrow{\cong} (2g+1) \langle
\omega, \omega \rangle \]
of line bundles on $\overline{\mathcal{I}}_g$. On the left hand side
one has a canonical non-zero rational section $\Lambda_g^{\otimes 3} \otimes
\delta^{\otimes - (2g+1)}$, and on the right hand side one has the canonical rational
section $q$. We claim that under $\mu$, these two rational sections are
identified, up to a sign. Indeed, the rational section $\Lambda_g^{\otimes 3} 
\otimes \delta^{\otimes -(2g+1)}$ restricts to the global section 
$\Lambda_g^{\otimes 3}$ of $3(8g+4)\lambda$ over $\mathcal{I}_g$ which is
nowhere vanishing. By the second half of 
(i) we have that $q$ is nowhere vanishing 
as well, and hence, over $\mathcal{I}_g$, the image of   
$\Lambda_g^{\otimes 3} \otimes \delta^{\otimes -(2g+1)}$ under $\mu$ differs
from $q$ by an invertible regular
function. It is stated in \cite[Lemma 7.3]{dj1} that such a function is either
$+1$ or $-1$. The claim follows.

The first half of (i) and statement (ii) follow from this 
fact that $\Lambda_g^{\otimes 3} \otimes
\delta^{\otimes -(2g+1)}$ and $q$ are identified, up to sign. As to the first
half of statement (i), 
from the Cornalba-Harris equality (\ref{ch}) we see that, 
if $S$ is normal, the rational section 
$\Lambda_g^{\otimes 3} \otimes \delta^{\otimes -(2g+1)}$ of $3(8g+4)\lambda-(2g+1)\delta$
is in fact global, as it gives rise to an
effective Cartier divisor on $S$. It follows that $q$ is global too.

As to (ii), first take $K$ to be a complete discrete valuation field. 
We may assume that
$X$ has semistable reduction over $K$. Let $R$ be the valuation ring of $K$ and
let $\mathcal{X} \to \mathrm{Spec}(R)$ be the regular minimal model of $X$. As
above we denote by $\nu(\cdot)$ order of vanishing along the closed point of
$\mathrm{Spec}(R)$. We have:
\[ \nu(q) = 3 \nu(\Lambda_g) - (2g+1) \nu(\delta) = 3 d(X)- (2g+1)\delta(X) \, .
\]
On the other hand, by \cite[Theorem 4.4]{zh1} we have:
\[ \nu(q) = -\log |q|_a + (2g+1) \varepsilon(X) \, , \]
and indeed the formula:
\[ -\log|q|_a = (2g-2) \chi(X) \]
drops out.

Next let $K$ be equal to $\mathbb{R}$ or $\mathbb{C}$. 
By \cite[Th\'eor\`eme 2.2]{mb} 
the isomorphism $\mu$, restricted to $\mathrm{Spec} (\bar{K})$, has 
admissible norm 
$\mathrm{e}^{(2g+1)\delta(X)}$. It follows that:
\[ - \log|q|_a = -3\log \| \Lambda_g \|(X) - (2g+1)\delta(X) = (2g-2) \chi(X) \]
in this case as well.

The formula in (iii) follows directly from the definition of $q$. Fix an index
$i$. Under the assumptions of (iii) we have:
\[ -\log |q|_a = -\log | (2^{4g} d_{ik} )^{g-1} \cdot q_{ik}^{\otimes 
2g+1} |_a  \]
for all $k \neq i$, hence, using the identity $\prod_{k \neq i } d_{ik} =1$ once more:
\[ - \log|q|_a = -4g(g-1) \log|2| - \sum_{k \neq i} \log |q_{ik}|_a \, . \]
We find:
\[ - \log |q|_a = -4g(g-1) \left( \log|2| + \sum_{k \neq i} (w_i,w_k)_a \right)
\]
as required.
\end{proof}
\begin{remark} Let $K$ be a number field or the function field of a curve over a
field which is not of characteristic~$2$. Let $X$ be a hyperelliptic curve over
$K$ of genus $g \geq 2$. Assume that $X$ has semistable reduction over $K$, and
that all Weierstrass points $w_1,\ldots,w_{2g+2}$ of $X$ are rational over $K$. Let
$(\omega,\omega)_a$ be the admissible self-intersection of the relative
dualising sheaf of $X$, and fix an index $i$. From (iii) of the above theorem 
we obtain, by summing over all places of $K$, that:
\[ (2g+1) (\omega,\omega)_a = -4g(g-1) \sum_{k \neq i} (w_i,w_k)_a \, , \]
where $(w_i,w_k)_a$ is now the global admissible pairing of $w_i,w_k$ on $X$.
We note that from the fact that the global $(,)_a$ restricts to minus 
the N\'eron-Tate pairing on degree-$0$ divisors \cite[Section 5.4]{zh1} 
one may actually infer the stronger identity:
\[ (\omega,\omega)_a = -4g(g-1) (w_i,w_k)_a \]
for any two indices $i,k$. It is not true, however, that the \emph{local}
intersections $(w_i,w_k)_a$ are in general independent of the choice of $i,k$.
\end{remark}
\begin{remark} The combinatorial optimisation 
methods used in \cite{ya2} to prove effective
lower bounds for $\chi$ in the non-archimedean case are quite complicated. It
would be interesting to see whether the simple formula from Theorem B could be
used to obtain good lower bounds for $\chi$ in an easier fashion. Note that the
global section $q$ in the function field case in some sense ``explains'' the
strict positivity of $(\omega,\omega)_a$ for non-isotrivial hyperelliptic
fibrations.
\end{remark}
\begin{remark} Let $K$ be the
fraction field of a discrete valuation ring $R$ of residue characteristic not
equal to~$2$, and fix a separable algebraic closure $\bar{K}$ of $K$. 
The symmetric equations $\mathcal{C}_{ij}^{\zeta}$ of a
hyperelliptic curve $X$ over $K$
have good properties, for example:
\begin{itemize}
\item if $X$ has potentially good reduction over $R$ then the coefficients of
$\mathcal{C}_{ij}^\zeta$ in $\bar{K}$ are integral over $R$, and 
$\mathcal{C}_{ij}^\zeta$
has good reduction over $R[\mathcal{C}_{ij}^{\zeta}]$;
\item if the coefficients of $\mathcal{C}_{ij}^\zeta$ in $\bar{K}$ are integral over
$R$, then $\mathcal{C}_{ij}^\zeta$ is a minimal equation of $X$ over
$R[\mathcal{C}_{ij}^\zeta]$.
\end{itemize}
For details and proofs of these facts we refer to \cite[Section 2]{gu}.
\end{remark}

\section{Connection with Zhang's $\varphi$-invariant} \label{phi}

Let $K$ be a field which is either a number field or the function field of a
curve over a field. 
Let $X$ be a smooth projective geometrically connected curve of
genus $g \geq 2$ over $K$, and assume that $X$ has semistable reduction over $K$. 
In a recent paper \cite{zh2} S.-W. Zhang introduced, for each place $v$ of $K$, 
a real-valued invariant $\varphi(X_v)$ of $X \otimes K_v$, as follows:
\begin{itemize}
\item if $v$ is a non-archimedean place, then:
\[ \varphi(X_v) = -\frac{1}{4} \delta(X_v) + \frac{1}{4} \int_{R(X_v)}
g_v(x,x)((10g+2) \mu_v - \delta_{K_{X_v}}) \, , \]
where:
\begin{itemize}
\item $\delta(X_v)$ is the number of singular points on the special
fiber of $X$ at $v$,
\item $R(X_v)$ is the reduction graph of $X$ at $v$, 
\item $g_v$ is the
Green's function for the admissible metric $\mu_v$ on $R(X_v)$,   
\item $K_{X_v}$ is the
canonical divisor on $R(X_v)$. 
\end{itemize} 
In particular, $\varphi(X_v)=0$ if $X$ has good reduction at $v$;
\item if $v$ is an archimedean place, then:
\[ \varphi(X_v) = \sum_{\ell} \frac{2}{\lambda_\ell} \sum_{m,n=1}^g  \left| 
\int_{X(\bar{K}_v)} \phi_\ell \omega_m \bar{\omega}_n \right|^2 \, , \]
where $\phi_\ell$ are the non-constant 
normalised real eigenforms of the Arakelov 
Laplacian on the compact Riemann surface $X(\bar{K}_v)$, with eigenvalues 
$\lambda_\ell>0$, and $(\omega_1,\ldots,\omega_g)$ is an orthonormal basis for 
the hermitian inner product 
$(\omega,\eta) \mapsto \frac{i}{2} \int_{X(\bar{K}_v)} \omega \, \bar{\eta}$
on the space of holomorphic differentials.
\end{itemize}
The adelic invariant $\varphi$ connects the
admissible self-intersection of the relative dualising sheaf of $X$ and the
height of the canonical Gross-Schoen cycle on $X^3$.
\begin{thm} (Zhang \cite{zh2}) Let $(\omega,\omega)_a$ be the admissible
self-intersection of the relative dualising sheaf of $X$, and let $\langle
\Delta_\xi,\Delta_\xi \rangle$ be the height of the canonical Gross-Schoen cycle
on $X^3$. Then:
\[ (\omega,\omega)_a = \frac{ 2g-2}{2g+1} \left( \langle \Delta_\xi,\Delta_\xi
\rangle + \sum_v \varphi(X_v) \log Nv \right) \, , \]
where $v$ runs over the places of $K$. Here the $Nv$ are the local factors
associated to $v$ as discussed in Section~\ref{chi} above. If $X$ is
hyperelliptic, then $\langle \Delta_\xi, \Delta_\xi \rangle =0$, and hence 
the formula:
\begin{equation} \label{eqnphi} 
(\omega,\omega)_a = \frac{ 2g-2}{2g+1} \sum_v \varphi(X_v) \log Nv  
\end{equation}
holds.
\end{thm}
The invariant $\varphi$ arises in a natural way as an adelic 
intersection number on the self-product $X^2$ of $X$, cf. \cite[Theorem
2.5.1]{zh2}. Note the striking resemblance of $\varphi$ with $\chi$, for
hyperelliptic curves: like $\varphi$, the invariant $\chi$ vanishes at each
non-archimedean place of good reduction, and one has (\ref{eqnphi}) with 
$\chi$ instead of $\varphi$, see formula (\ref{formulachi}).

The canonical nature of both $\varphi$ and $\chi$ leads us to make the
following:
\begin{conj} \label{conj}
Assume that $X$ is hyperelliptic. Then
$\varphi(X_v)=\chi(X_v)$ for each place $v$ of $K$.
\end{conj}
This conjecture has the following interesting corollary. 
Assume that $K$ is a number field. Let $\pi \colon 
\mathcal{X} \to S$ be the regular minimal model of $X$ over 
the ring of integers of $K$. Then $\pi$ is a generically smooth semistable
hyperelliptic curve. By Theorem \ref{chithm}(i) we have $q$ as a canonical
global section of the bundle $(2g+1) \langle \omega, \omega \rangle$ on
$S$, and by (ii) of the same theorem   
$ (2g-2)\chi(X_v)=-\log|q|_a$ at
each place $v$ of $K$. The latter quantity is, as we have seen from the results of
Moriwaki and Yamaki, non-negative if $v$ is non-archimedean.

Consider now the case that $v$ is an archimedean place. It is
clear from the definition of $\varphi$ given above 
that $\varphi(X_v) \geq 0$. But in fact
one has $\varphi(X_v) >0$, as is explained in \cite[Remark 2.5.1]{zh2}. So, if
$\varphi(X_v)=\chi(X_v)$ would hold, it follows that $-\log|q|_a $ is 
\emph{positive} on the archimedean places. 
In other words, the section $q$ is a \emph{small} section of 
$(2g+1) \langle \omega, \omega \rangle$. The Bogomolov conjecture, i.e. the
statement that $(\omega,\omega)_a$ is strictly positive, would follow in a
conceptual way by taking the admissible Arakelov degree of $q$.\\

\noindent \emph{Note, not included in the published version.} It follows from \cite[Theorem 3.5]{ya3} (non-archimedean case) and \cite[Corollary 1.8]{dj2} (archimedean case) that Conjecture \ref{conj} is true.

\section*{Acknowledgments}
The author is supported by VENI grant 639.033.402
from the Netherlands Organisation for Scientific Research (NWO). He thanks the
Max Planck Institut f\"ur Mathematik in Bonn for its hospitality during a visit.

\vspace{.5cm}

\noindent Address of the author: \\  \\
Robin de Jong \\
Mathematical Institute \\
University of Leiden \\
PO Box 9512 \\
2300 RA Leiden \\
The Netherlands \\
Email: \verb+rdejong@math.leidenuniv.nl+


\begin{thebibliography}{99}

\bibitem{bo} S. Bosch, \emph{Formelle Standardmodelle hyperelliptischer 
Kurven}. Math. Ann. 251 (1980), no. 1, 19--42. 

\bibitem{bmmb} J.-B. Bost, J.-F. Mestre, L. Moret-Bailly, \emph{Sur le 
calcul explicite des ``classes de Chern'' des surfaces arithm\'etiques de 
genre $2$}. 
S\'eminaire sur les Pinceaux de Courbes Elliptiques (Paris, 1988).
Ast\'erisque No. 183 (1990), 69--105. 

\bibitem{fa} G. Faltings, \emph{Calculus on arithmetic surfaces}.  
Ann. of Math. 119  (1984),  no. 2, 387--424. 

\bibitem{gu} J. Gu\`ardia, \emph{Jacobian Nullwerte, periods and 
symmetric equations for hyperelliptic curves}.  Ann. Inst. Fourier 
(Grenoble)  57  (2007),  no. 4, 1253--1283. 

\bibitem{dj1} R. de Jong, \emph{Explicit Mumford isomorphism for hyperelliptic
curves}.  J. Pure Appl. Algebra  208  (2007),  no. 1, 1--14.

\bibitem{dj2} R. de Jong, \emph{Second variation of Zhang's $\lambda$-invariant on the moduli space of curves}. To appear in Amer. Jnl. Math.

\bibitem{ka} I. Kausz, \emph{A discriminant and an upper bound for 
$\omega^2$ for hyperelliptic arithmetic surfaces}.  
Compositio Math.  115  (1999),  no. 1, 37--69.

\bibitem{ll} O.A. Laudal and K. L\o nsted,
\emph{Deformations of curves. I. Moduli for hyperelliptic curves.}  
Algebraic geometry (Proc. Sympos., Univ. Troms\o, Troms\o, 1977) 
150--167, Lecture Notes in Mathematics 687, Springer, Berlin, 1978.

\bibitem{lo} P. Lockhart, \emph{On the discriminant of a hyperelliptic 
curve}. Trans. Amer. Math. Soc. 342 (1994), no. 2, 729--752. 

\bibitem{lk} K. L\o nsted, S.L. Kleiman, 
\emph{Basics on families of hyperelliptic curves}.
Compositio Math. 38 (1979), no. 1, 83--111. 

\bibitem{ma} S. Maugeais, \emph{Rel\`evement des rev\^etements $p$-cycliques 
des courbes rationnelles semi-stables}.  
Math. Ann. 327 (2003), no. 2, 365--393. 

\bibitem{mb} L. Moret-Bailly, \emph{La formule de Noether pour les surfaces 
arithm\'etiques}.  Invent. Math. 98 (1989), no. 3, 491--498. 

\bibitem{mo} A. Moriwaki, 
\emph{Bogomolov conjecture for curves of genus $2$ over function fields}.  
J. Math. Kyoto Univ.  36  (1996),  no. 4, 687--695.

\bibitem{sga1} \emph{S\'eminaire de g\'eom\'etrie alg\'ebrique du Bois Marie
1960-61 I: Rev\^etements \'etales et groupe fondamental}. Dirig\'e par A.
Grothendieck. Documents Math\'ematiques 3, Soci\'et\'e Math\'ematique de France
2003.

\bibitem{ya1} K. Yamaki, \emph{Cornalba-Harris equality for semistable 
hyperelliptic curves in positive characteristic}. 
 Asian J. Math.  8  (2004),  no. 3, 409--426. 

\bibitem{ya2} K. Yamaki, \emph{Effective calculation of the geometric 
height and the Bogomolov conjecture for hyperelliptic curves over 
function fields}.  J. Math. Kyoto Univ.  48  (2008),  no. 2, 401--443. 

\bibitem{ya3} K. Yamaki, \emph{Graph invariants and the positivity of the height of the
Gross-Schoen cycle for some curves}. manuscripta math. 131 (2010), 149--177.

\bibitem{zh1} S.-W. Zhang, \emph{Admissible pairing on a curve}.
Invent. Math.  112  (1993),  no. 1, 171--193.

\bibitem{zh2} S.-W. Zhang, \emph{Gross-Schoen cycles and dualising sheaves}. 
Invent. Math. 179 (2010), 1--73.

\end{thebibliography}
\end{document}